\documentclass[preprint,12pt]{elsarticle}

\usepackage{amssymb, amsmath}
\usepackage{color} 

\newproof{pf}{Proof}
\newtheorem{definition}{Definition}
\newtheorem{theorem}{Theorem}

\begin{document}

\begin{frontmatter}



\title{Simultaneous Approximation of a Multivariate Function and its Derivatives by Multilinear Splines}


\author[Anderson]{Ryan~Anderson}

\author[Babenko]{Yuliya~Babenko\corref{cor1}}

\author[Leskevych]{Tetiana~Leskevych}

\cortext[cor1]{Corresponding author}

\address[Anderson]{Department of Mathematics and Statistics, Kennesaw State University, 1000 Chastain Road, \#1601, Kennesaw, GA 30144-5591, USA}
\address[Babenko]{Department of Mathematics and Statistics, Kennesaw State University, 1000 Chastain Road, \#1601, Kennesaw, GA 30144-5591, USA}
\address[Leskevych]{Department of Mathematics and Mechanics, Dnepropetrovsk National University, Gagarina pr., 72, Dnepropetrovsk, 49010, UKRAINE}

\begin{abstract}
In this paper we consider the approximation of a function by its interpolating multilinear spline and the approximation of its derivatives by the derivatives of the corresponding spline.
We derive formulas for the uniform approximation error on classes of functions with moduli of continuity bounded above by certain majorants.
\end{abstract}


\end{frontmatter}



\section{Basic definitions and notation}
\indent  Let ${\bf x}=(x_1,x_2,...,x_n)$ be a point in Euclidean space $\mathbb{R} ^n$. By $C_D$ we denote the class of functions $f({\bf x})=f(x_1,x_2,...,x_n)$ that are continuous on the domain $D:=[0,1]^n\subset \mathbb{R} ^n$.


We consider a vector ${\bf r}\in \{0,1\}^n$, i.e. a vector having $n$ components each being either $0$ or $1$.
Let
$C_D^{{\bf r}}$ be the class of functions, $f({\bf x})\in C_D$, with continuous derivatives
\[
f^{({\bf t})}({\bf x})=\frac{\partial^{\sum_{i=1}^{n}t_{i}}f } {\partial x^{t_1}...\partial x^{t_n}}{ ({\bf x})},
\]
where ${\bf t}\in \{0,1\}^n$ and $t_i\leq r_i$ for each $i=1,...,n$. We define $f^{({\bf 0})}({\bf x}):=f({\bf x})$ and $C_D^{{\bf 0}}:=C_D$.\\
\indent For any function $f({\bf x})\in C_D$, consistent with literature, we denote the uniform
norm as
\[
\|f \|_C:=\max \{ |f({\bf x})|:{\bf x}\in D \}.
\]
%
\indent 
The next two definitions introduce two types of moduli of continuity of a given function $f$, both
 characterizing the smoothness of the original function $f$.

\begin{definition}\label{defomsmmul}
If the function $f({\bf x})$ is bounded for $x_i\in [a_i,b_i], i=1,...,n$, then its total modulus of continuity, $\omega(f;\boldsymbol\tau)$, is defined as follows
\begin{eqnarray*}
\omega(f;\boldsymbol\tau)& := & \omega(f;\bf{a},\bf{b};\boldsymbol\tau)\\
& = & \sup\left\{|f({\bf x})-f({\bf y})|:\;\;|x_i-y_i|\leq \tau_i;\;x_i,y_i\in [a_i,b_i] \right\},
\end{eqnarray*}
where $0\leq \tau_i \leq b_i-a_i$, for $i=1,...,n$ and $\boldsymbol\tau:=(\tau_1,...,\tau_n)$, $\boldsymbol a:=(a_1,...,a_n)$, $\boldsymbol b:=(b_1,...,b_n)$.\\
\end{definition}
In addition, we consider the following $l_p$  distances, $1\le p<\infty$, between points ${\bf x,y}\in D\subset \mathbb{R}^n$,
$$
\|{\bf x}-{\bf y}\|_p:=\sqrt[p]{\sum_{i=1}^{n}|x_i-y_i|^p}.
$$
{
\begin{definition}\label{defomsmone}
For the function $f({\bf x})\in C_D$ and for given $p$, $1\le p<\infty$,  we define the modulus of continuity of function $f$ with respect to $p$ to be
$$
\omega_{p}\left(f;\gamma  \right)
  := \sup\left\{|f({\bf x})-f({\bf y})|:\; {\bf x,y}\in D,\;\; \|{\bf x}-{\bf y}\|_p\leq\gamma \right\},\;0\leq \gamma \leq d_p,
$$
where $d_p:=\max\{ \|{\bf x}-{\bf y}\|_p:\;\; {\bf x,y}\in D\subset \mathbb{R}^n \}$.
\end{definition}
}
We point out that
$$
{d}_{p}=\sqrt[p]{n}.
$$
Note that the moduli of continuity of all suitable functions have some common properties. We call all functions (univariate or multivariate) with these properties {\it functions of the moduli of continuity type} and use them to define classes of smoothness of functions.
%
\begin{definition}\label{defombig}
Function $\Omega(\boldsymbol\tau)$ is called a function of modulus of continuity type, or MC-type function (for short),
if the following properties hold for any vectors $\boldsymbol \gamma, \boldsymbol\tau \in \mathbb{R}^n_{+}:=
\{{\bf x}\in\mathbb{R}^n:x_i\geq0,\;i=1,...,n\}$:
\begin{enumerate}
\item $\Omega({\bf 0})=0$.
\item $\Omega(\boldsymbol\tau)$ is non-decreasing.
\item $\Omega(\boldsymbol\tau+\boldsymbol \gamma) \leq \Omega(\boldsymbol\tau)+\Omega(\boldsymbol \gamma)$, that is $\Omega(\boldsymbol\tau)$ is subadditive.
\item $\Omega(\boldsymbol\tau)$ is continuous for all $\tau_i
$, $i=1,...,n$.
\end{enumerate}
\end{definition}

The following two definitions present the classes of functions that we will be working with in this paper.

\begin{definition}\label{defclassmul}
Given an arbitrary $n$-variate MC-type function $\Omega({\boldsymbol \tau})$,  we define the class
$C_D^{\bf r}(\Omega)$, where ${\bf r}\in \{0,1\}^n$ with
$C_D^{{\bf 0}}(\Omega):=C_D(\Omega)$, to be the class of functions $f({\bf x})\in C_D^{\bf r}$, such that the total modulus of continuity of their derivatives of order ${\bf r}$ are bounded from above by the given $\Omega({\boldsymbol \tau})$
\begin{eqnarray}\label{omega}
\omega \left(f^{({\bf r})}; {\bf \boldsymbol\tau} \right)&
\leq & \Omega({\boldsymbol \tau}),\qquad 0\leq \tau_i\leq 1,\;\; i=1,...,n.
\end{eqnarray}
\end{definition}

\begin{definition}\label{defclassone}
 Given an arbitrary univariate MC-type function $\Omega({\gamma})$, we define the class
$C_{D,p}^{\bf r}(\Omega)$, $1\leq p \leq \infty$, where ${\bf r}\in \{0,1\}^n$ with
$C_{D,p}^{{\bf 0}}(\Omega):=C_{D,p}(\Omega)$, to be the class of functions $f({\bf x})\in C_D^{\bf r}$, such that their moduli of continuity are bounded from above by the given $\Omega({\gamma})$
\begin{equation}\label{Omega}
\omega_{p}\left(f^{({\bf r})};\gamma  \right)
\leq \Omega(\gamma),\qquad 0\leq \gamma \leq \sqrt[p]{n}.
\end{equation}
\end{definition}



\section{Construction of the interpolating spline}
\indent In order to construct the interpolating spline for the given function $f\in C_D$, which will be the main approximation tool for $f$ (and its derivatives will be used to approximate the derivatives of $f$), we fix a vector ${\bf m}=(m_1,...,m_n)\in\mathbb{N}^n$ and for each $i=1,...,n$ we first define the univariate grid of nodes as follows
$$
{\cal D}_{m_i}=\left \{ 0, \frac{1}{m_i},..., \frac{m_i-1}{m_i},1\right \}.
$$
With the help of the standard Cartesian product we define the $n$-variate grid as
$$
{\cal D}_{{\bf m}}={\cal D}_{m_1}\times ...\times {\cal D}_{m_n}.
$$
Once the grid is constructed, each point on the grid is defined by a vector
${\bf j}=(j_1,\dots ,j_n)$ where $j_i\in \{0,\dots, m_i\}$,
as follows
$$
{\bf x}^{\bf j}:=(x_1^{j_1},\dots, x_n^{j_n}):=\left (\frac{j_1}{m_1},\dots , \frac{j_n}{m_n}\right )\in {\cal D}_{{\bf m}}.
$$

%


Having defined the grid,
we next define the interpolating spline that will be used to approximate the given function and whose derivatives will be used to simultaneously approximate  the derivatives of the function.
\begin{definition}
For the grid of nodes, ${\cal D}_{{\bf m}}$, and a given function $f({\bf x})\in C_D$, we define the multilinear ($n$-linear) interpolating spline, $S_{\bf m}(f;{\bf x})$, to satisfy the following conditions:
\begin{enumerate}
\item On every block $D_{{\bf j}}:=\prod_{i=1}^{n}[x_i^{j_{i}},x_i^{j_{i}+1}]$, where ${\bf j}=(j_1,...,j_n)$, $j_i=0,...,m_i-1$, $S_{\bf m}(f;{\bf x})$ is an algebraic polynomial of first degree in $x_i$ for $i=1,...,n$.
\item $S_{\bf m}(f;{\bf x}^{\bf j})=f({\bf x}^{\bf j})$ for ${\bf j}=(j_1,...,j_n)$, $j_i=0,...,m_i$ and $i=1,...,n$. In other words, $S_{\bf m}(f;{\bf x})$ {\it interpolates} $f({\bf x})$ at the nodes ${\cal D}_{{\bf m}}$.\\
\end{enumerate}
\end{definition}

\indent Note that for ${\bf x}\in D_{{\bf j}}$, $j_i=0,...,m_i-1$ and $i=1,...,n$, the following holds
\begin{equation}\label{spline}
S_{\bf m}(f;{\bf x})=\sum_{l_1=0}^{1}...\sum_{l_n=0 }^{1}\left[ f({\bf x}^{{\bf j}+{\bf l}})\left(\prod_{i=1}^{n}H_{l_i,j_{i}}(x_i)\right) \right],
\end{equation}
where ${\bf l}=(l_1,...,l_n)$ and
\begin{equation}\label{prop H}
H_{0,j_i}(x_i):=m_i(x_i^{j_{i}+1}-x_i),\qquad \sum_{l_{i}=0}^1 H_{l_i,j_i}(x_i)\equiv 1.
\end{equation}

\indent When taking the partial derivatives of $S_{\bf m}(f;{\bf x})$, there may be discontinuities at the points of the following set
\begin{equation*}
A:=\{{\bf x}: \;\; x_i^{j_{i}} \mbox{ is a component of } {\bf x} \mbox{ for some } i\; 
 {\mbox {and some}} \; j_i\in\{0,...,m_i\}\}.
\end{equation*}
\indent The discontinuities of the partial derivatives of $S_{\bf m}(f;{\bf x})$ may exist because $S_{\bf m}(f;{\bf x})$ is a piecewise linear polynomial with respect to $x_i$, for $i=1,...,n$.  With these discontinuities all partial derivatives need to be defined carefully.  In order to do so, we define the set $M$
of indices as follows:
\begin{eqnarray*}
M & := & \{i:\;\;\; \mbox{The derivative of } S_{\bf m}(f;{\bf x}) \mbox{ is taken with respect to } x_i\}.
\end{eqnarray*}
\indent We use the notation $|M|$ to denote the number of elements (cardinality) in set $M$.

\indent Next, we introduce the functions $F_{M}({\bf x})$ and $F_{M'}({\bf x})$ - which will later be used to define the partial derivatives of $S_{\bf m}(f;{\bf x})$ - as follows:
\begin{eqnarray*}
F_{M}({\bf x})& := & \sum_{l_1=0}^{1}...\sum_{l_n=0 }^{1}(-1)^{\sum_{i\in M}l_i}\left(\prod_{i\notin M}H_{l_i,j_i}(x_i) \right) f({\bf x}^{{\bf j}+{\bf l}}),\\
F_{M'}({\bf x})& := & \sum_{l_1=0}^{1}...\sum_{l_n=0 }^{1}(-1)^{\sum_{i\in M}l_i+1}\left(\prod_{i\notin M}H_{l_i,j_i}(x_i) \right) f({\bf x}^{{\bf j}+{\bf l}})\nonumber.
\end{eqnarray*}
{\bf Remark.} $F_{M}({\bf x})$ is used if $|M|$ is even, and $F_{M'}({\bf x})$ is used if $|M|$ is odd. When stating the results,
we will use  {M such that $|M|$ is even, but the results also hold for M in the case when $|M|$ is odd.}

Let ${\bf b}\in \{0,1\}^n$, where
\begin{equation*}
b_i:=\left\{ \begin{array}{lr}
1, & \mbox{if } i\in M\\
0, & \mbox{if } i\notin M
\end{array}   \right. ,
\end{equation*}
then
\begin{equation}\label{divspline}
S_{\bf m}^{({\bf b})}(f;{\bf x})=\left(\prod_{i\in M}m_i\right)F_{M}({\bf x}),
\end{equation}
for ${\bf x}\in D_{{\bf j}}^{'}:=\prod_{i=1}^{n}I(x_i)$, where
\begin{eqnarray*}
I(x_i):= \left\{
\begin{array}{ll}
\left[ x_i^{j_{i}},x_i^{j_{i}+1}\right), & \mbox{if } j_i=0,...,m_i-2\\
\left[x_i^{j_{i}},x_i^{j_{i}+1}\right],&\mbox{if } j_i=m_i-1
\end{array}. \right.
\end{eqnarray*}

Finally, we introduce the errors of approximation for a given function (Definition \ref{deferr}) and error of approximation on a class of functions (Definition~\ref{deferrclass}).
\begin{definition}\label{deferr}
For {given ${\bf r}\in \{0,1\}^n$ and} a function $f\in C_{D}^{{\bf r}}$,
we denote the {\it error of approximation of a function} $f^{({\bf r})}({\bf x})$ (or its derivative) by the interpolating spline $S_{\bf m}^{({\bf r})}(f;{\bf x})$ (or its derivative, respectively) constructed above to be
\begin{equation}\label{error}
{\cal E}_{{\bf m}}^{{\bf r}}(f;{\bf x}):=\left|f^{({\bf r})}({\bf x})-S_{\bf m}^{({\bf r})}(f;{\bf x})\right|,\;\;\; {\bf x}\in D
\end{equation}
with ${\cal E}_{{\bf m}}^{{\bf 0}}(f;{\bf x}):={\cal E}_{{\bf m}}(f;{\bf x})$.\\
\end{definition}
\begin{definition}\label{deferrclass}
For a {given ${\bf r}\in \{0,1\}^n$} and for any class ${\cal M}\in C_{D}^{{\bf r}}$, we denote the {\it error of
approximation on the class ${\cal M}$ by splines that interpolate at the nodes ${\cal D}_{{\bf m}}$}
  to be
\begin{equation}\label{maxerror}
E_{{\bf m}}^{{\bf r}}({\cal M}):=\sup \left\{\|{\cal E}_{{\bf m}}^{{\bf r}}(f)\|_C\; :\;\; f\in{\cal M}  \right\}
\end{equation}
with $E_{{\bf m}}^{{\bf 0}}({\cal M}):=E_{{\bf m}}({\cal M})$.
\end{definition}

In this paper, we present the explicit formulas for the uniform error of approximation of multivariate functions from some classes of smoothness by
multilinear interpolating splines as well as the error of approximation of the derivatives of functions from the
considered class by the derivatives of the {corresponding} splines. An analogous univariate result for approximating the function from the same class is contained in  ~\cite{Malozemov66}, and the result for approximating the derivatives is contained in ~\cite{Malozemov67}.
In the case of bivariate functions from the class $C_{D}^{\bf r}(\Omega)$, earlier known results for such functions are in the paper of Storchai ~\cite{Storchai75} and for the derivatives are in works of Vakarchuk and Shabozov ~\cite{Vakarchuk90, Shabozov}. For the class $C_{D,p}^{\bf r}(\Omega)$ the known results are only for the cases $p=1$ ~\cite{Vakarchuk05} and $p=2$ ~\cite{Storchai72} for functions and
 ~\cite{Vakarchuk05,Shabozov} for derivatives, respectively. We have extended their results to the case of arbitrary dimension and arbitrary $1\leq p\leq 3$.

%

\section{The error of approximation on classes $C_D(\Omega)$ and $C_{D,p}(\Omega)$ }
\indent In this section, we estimate the error of approximation by interpolating splines
on the classes $C_D(\Omega)$ and $C_{D,p}(\Omega)$ defined in Definitions \ref{defclassmul} and \ref{defclassone}.

\begin{theorem}\label{thm1}
Let $\Omega(\boldsymbol\tau)$ be an arbitrary concave (in each variable) MC-type function. Then for ${\bf m}\in \mathbb{N}^n$ with $m_i\geq2$ for $i=1,...,n$, the error on the class $C_D(\Omega)$ is
\[
E_{{\bf m}}\left(C_D(\Omega) \right)=\Omega \left(\frac{1}{2{\bf m}}
\right).
\]\end{theorem}
{\bf Proof}\\
Let an arbitrary function $f\in C_D(\Omega)$ be given.
 Without loss of generality we consider ${\bf x}\in D_{{\bf j}}=\prod_{i=1}^{n}[x_i^{j_{i}},x_i^{j_{i}+1}]$ for some ${\bf j}=(j_1,...,j_n)$, ${j_i\in\{0,...,m_i-1\}}$. By Definition \ref{deferr}, using connection \eqref{prop H}, we have
\begin{eqnarray*}
{\cal E}_{{\bf m}}(f;{\bf x})& = &
|f({\bf x})-S_{\bf m}(f;{\bf x})|\\
& = & \left|\sum_{l_1=0}^1 ... \sum_{l_n=0}^1 \left( \prod_{i=1}^{n} H_{l_i,j_i}(x_i) \right) \left(f({\bf x})-f({\bf x}^{{\bf j}+{\bf l}})\right)\right|,
\end{eqnarray*}
where ${\bf x}^{{\bf j}+{\bf l}}=(x_1^{j_{1}+l_{1}},...,x_n^{j_{n}+l_{n}})$ and 
${\bf l}=(l_1,...,l_n)\in\{0,1\}^n$ for ${i=1,...,n}$.\\
\indent Using the triangle inequality, Definition \ref{defomsmmul}, and Definition \ref{defclassmul}, we obtain
\begin{eqnarray*}
|{\cal E}_{{\bf m}}(f;{\bf x})| & \leq &  \sum_{l_1=0}^1 ... \sum_{l_n=0}^1 \left( \prod_{i=1}^{n} H_{l_i,j_i}(x_i) \right) \left|f({\bf x})-f({\bf x}^{{\bf j}+{\bf l}})\right| \\
& \leq & \sum_{l_1=0}^1 ... \sum_{l_n=0}^1 \left( \prod_{i=1}^{n} H_{l_i,j_i}(x_i) \right) \omega(f;|{\bf x}-{\bf x}^{{\bf j}+{\bf l}}|)\\
& \leq & \sum_{l_1=0}^1 ... \sum_{l_n=0}^1 \left( \prod_{i=1}^{n} H_{l_i,j_i}(x_i) \right) \Omega(|{\bf x}-{\bf x}^{{\bf j}+{\bf l}}|),
\end{eqnarray*}
where $|{\bf x}-{\bf x}^{{\bf j}+{\bf l}}|=(|x_1-x_1^{j_{1}+l_{1}}|,...,|x_n-x_n^{j_{n}+l_{n}}|)$.\\
\indent Since ${\bf x}\in D_{{\bf j}}$ and function $\Omega(\boldsymbol\tau)$ is concave in each variable, we have
\begin{eqnarray*}
|{\cal E}_{{\bf m}}(f;{\bf x})|
& \leq & \Omega\left( \lambda(x_1),...,\lambda(x_n) \right),
\end{eqnarray*}
where
\begin{equation}\label{lambda}
\lambda(x_i):=H_{0,j_i}(x_i)(x_i-x_i^{j_{i}})+H_{1,j_i}(x_i)(x_i^{j_{i}+1}-x_i),\;\;\; i=1,...,n.
\end{equation}

\indent Since the function of MC-type is non-decreasing, we need to find
$\max\{\lambda(x_i): x_i\in [x_i^{j_{i}},x_i^{j_{i}+1}]  \}$ for $i=1,...,n$.
For convenience, we use the substitution $x_i=t$:
\begin{eqnarray*}
\lambda(t)&=& H_{0,j_i}(t)(t-x_i^{j_{i}})+H_{1,j_i}(t)(x_i^{j_{i}+1}-t)\\
& = & (m_ix_i^{j_{i}+1}-m_it)(t-x_i^{j_{i}})+(1-m_ix_i^{j_{i}+1}+m_it)(x_i^{j_{i}+1}-t)\\
& = & 3m_ix_i^{j_{i}+1}t-2m_it^2+m_ix_i^{j_{i}}t-m_ix_i^{j_{i}}x_i^{j_{i}+1}-m_i\left(x_i^{j_{i}+1}\right)^2+x_i^{j_{i}+1}-t.
\end{eqnarray*}
\indent Thus, taking the derivative of $\lambda(t)$ with respect to $t$,
setting $\lambda'(t)=0$, and solving for $t$ yields
$$
t= \frac{x_i^{j_{i}}+x_i^{j_{i}+1}}{2}, \qquad \hbox{and hence} \qquad \lambda\left(\frac{x_i^{j_{i}}+x_i^{j_{i}+1}}{2}\right)=\frac{1}{2m_i}.
$$
In order to show that $\frac{1}{2m_i}$ is a maximum, we take second derivative 
\[
\lambda''(t)=-4m_i.
\]
Since $m_i\geq 2$, we have $\lambda''(t)<0$. Hence, $\frac{1}{2m_i}$ is a maximum and we obtain
\begin{equation}\label{lambdamax}
\max\{ \lambda(x_i):x_i\in [x_i^{j_{i}},x_i^{j_{i}+1}] \}=\lambda\left( \frac{x_i^{j_{i}}+x_i^{j_{i}+1}}{2}\right)=\frac{1}{2m_i},\;\;\;i=1,...,n.
\end{equation}

\indent As $\Omega(\boldsymbol\tau)$ is non-decreasing, we conclude
\begin{equation*}
|{\cal{E}}_{{\bf m}}\left(f;x \right)|\leq \Omega \left(
\frac{1}{2{\bf m}}
\right).
\end{equation*}
Since last inequality holds for any $f\in C_D(\Omega)$, we have
\begin{equation}\label{preresult1}
E_{{\bf m}}\left(C_D(\Omega) \right)\leq \Omega \left(
\frac{1}{2{\bf m}}
\right).
\end{equation}
\indent Next, we present a particular function {from $C_D(\Omega)$} (called an {\it extremal function}) for which (\ref{preresult1}) occurs with equality. We define the extremal function $f_0$ as follows
\begin{eqnarray*}
f_0({\bf x})&:=&\Omega\left( \min\{x_1-x_1^{j_1},x_1^{j_1+1}-x_1\},...,\min\{x_n-x_n^{j_n},x_n^{j_n+1}-x_n\}\right),\\
\end{eqnarray*}
where ${\bf x}=(x_1,...,x_n) \in D_{{\bf j}}$ {for fixed ${\bf j}$} and $f_0({\bf x})=0$ for ${\bf x}\notin D_{{\bf j}}$. Then we have $f_0\in C_D(\Omega)$, and $f_0({\bf x}^{\bf j})=0$ for ${\bf j}=(j_1,...,j_n)$, $j_i=0,...,m_i$ and $i=1,...,n$. In addition, we have
\[
\|{\cal E}_{{\bf m}}(f_0) \| _C = \|f_0\|_C,
\]
since $S_{\bf m}(f_0;{\bf x})$ is linear on each partition element and interpolates $f_0({\bf x})$ at the nodes of ${\cal D}_{{\bf m}}$. Consequently, the following estimates hold:
\begin{eqnarray}\label{EEE}
E_{{\bf m}}\left(C_D(\Omega) \right) & \geq &\|{\cal E}_{{\bf m}}(f_0)\| _C
 =  \|f_0\|_C  
 \geq  f_0\left( 
\frac{1}{2{\bf m}}
\right)
 =  \Omega \left( 
\frac{1}{2{\bf m}}
\right). 
\end{eqnarray}
\indent Combining (\ref{preresult1}) and (\ref{EEE}), we obtain
\[
E_{{\bf m}}\left(C_D(\Omega) \right)=\Omega \left( 
\frac{1}{2{\bf m}}
\right). \;\;\;\; \square
\]

{
Theorem \ref{thm2} gives the error $E_{{\bf m}}(C_{D,p}(\Omega))$ for
$1 \leqslant p \leqslant 3$, where $C_{D,p}(\Omega)$ is defined in (\ref{Omega}).
\begin{theorem}\label{thm2}
Let $\Omega(\gamma)$ be an arbitrary concave, MC-type, univariate function. Then for ${\bf m}\in \mathbb{N}^n$ with $ m_i\geq 2$ for all $i=0,...,n$, the error on the class $C_{D,p}(\Omega)$ for $ 1 \leqslant p \leqslant 3$ is
\[
E_{{\bf m}}\left(C_{D,p}(\Omega) \right)=\Omega \left(\frac{1}{2}\sqrt[p]{\sum\limits_{i=1}^n\frac{1}{m_i^p}}\right).
\]
\end{theorem}
{\bf Proof}

\indent For any $\bf x$ from an arbitrary $D_{{\bf j}}:=\prod_{i=1}^{n}[x_i^{j_{i}},x_i^{j_{i}+1}]$ and for an arbitrary $f\in C_{D,p}(\Omega)$ by Definition \ref{deferr}, we have
\begin{eqnarray*}
{\cal E}_{{\bf m}}(f;{\bf x})& = &\left|f({\bf x})-S_{\bf m}(f;{\bf x})\right|\\
& \le & \sum_{l_1=0}^1 ... \sum_{l_n=0}^1 \left( \prod_{i=1}^{n} H_{l_i,j_i}(x_i) \right) \left|f({\bf x})-f({\bf x}^{{\bf j}+{\bf l}})\right|,
\end{eqnarray*}
where ${\bf x}^{{\bf j}+{\bf l}}=(x_1^{j_{1}+l_{1}},...,x_n^{j_{n}+l_{n}})$ and $j_i=0,...,m_i-1$, ${\bf l}=(l_1,...,l_n)\in\{0,1\}^n$ for ${i=1,...,n}$.\\
\indent Using Definition \ref{defomsmone}, we have
\begin{equation*}
{\cal E}_{{\bf m}}(f;{\bf x})  \leq  \sum_{l_1=0}^1 ... \sum_{l_n=0}^1 \left( \prod_{i=1}^{n} H_{l_i,j_i}(x_i) \right)\omega_{p}\left(\sqrt[p]{\sum_{i=1}^{n}|x_i-x_i^{j_{i}+l_{i}}|^p}  \right).
\end{equation*}
\indent Using Definition \ref{defclassone} for the class $C_{D,p}(\Omega)$, relation (\ref{prop H}), and the fact that $\Omega(\gamma)$ and function $t^{\frac{1}{p}}$, $ 1 \leqslant p \leqslant 3$, are both concave functions, we obtain

\begin{eqnarray*}
{\cal E}_{{\bf m}}(f;{\bf x})& \leq & \sum_{l_1=0}^1 ... \sum_{l_n=0}^1 \left( \prod_{i=1}^{n} H_{l_i,j_i}(x_i) \right)\Omega\left(\sqrt[p]{\sum_{i=1}^{n}|x_i-x_i^{j_{i}+l_{i}}|^p}  \right)\\
& \leq & \Omega \left( \sqrt[p]{\sum_{i=1}^{n}\left[ \sum_{l_1=0}^1 ... \sum_{l_n=0}^1 \left( \prod_{i=1}^{n} H_{l_i,j_i}(x_i) \right)|x_i-x_i^{j_{i}+l_{i}}|^{p}\right]}  \right)\\
\end{eqnarray*}
\vspace{-6,5ex}
\begin{eqnarray*}
& = & \Omega \left( \sqrt[p]{\sum_{i=1}^{n}\left[ \left( \sum_{l_i=0}^{1}H_{l_i,j_i}(x_i)|x_i-x_i^{j_{i}+l_{i}}|^p\right) \left(\prod_{ k=1;\;\;k\neq i}^{n}\left[\sum_{l_k=0}^{1}H_{l_k,j_k}(x_k) \right]  \right)\right]}\right)  \\
& = & \Omega \left( \sqrt[p]{\sum_{i=1}^{n}\left[ \sum_{l_i=0}^{1}H_{l_i,j_i}(x_i)|x_i-x_i^{j_{i}+l_{i}}|^p \right]}\right).
\end{eqnarray*}
\indent For $i=1,...,n$, we denote
\begin{eqnarray}\label{alpha}
\alpha(x_i)&:=& \sum_{l_i=0}^{1}H_{l_i,j_i}(x_i)|x_i-x_i^{j_{i}+l_{i}}|^p\nonumber\\
&=&H_{0,j_i}(x_i)(x_i-x_i^{j_{i}})^p+H_{1,j_i}(x_i)(x_i^{j_{i}+1}-x_i)^p.
\end{eqnarray}
We have
\begin{eqnarray*}
\alpha(x_i):=H_{0,l_i}(x_i)(x_i-x_i^{j_{i}})^p+H_{1,l_i}(x_i)(x_i^{j_{i}+1}-x_i)^p\\
=m_i(x_i^{j_{i}+1}-x_i)(x_i-x_i^{j_{i}})^p+m_i\left(\frac 1 {m_i}-x_i^{j_{i}+1}+x_i\right)(x_i^{j_{i}+1}-x_i)^p\\
=m_i(x_i^{j_{i}+1}-x_i)(x_i-x_i^{j_{i}})^p+m_i(x_i-x_i^{j_{i}})(x_i^{j_{i}+1}-x_i)^p\\
=m_i(x_i^{j_{i}+1}-x_i)^{p+1}\left(\left(\frac{x_i-x_i^{j_{i}}}{x_i^{j_{i}+1}-x_i}\right)^p+\frac{x_i-x_i^{j_{i}}}{x_i^{j_{i}+1}-x_i}\right)
\end{eqnarray*}
Using the following inequality \cite[p. 334]{Korneychuk}
$$
{ 2}^{{ p}} \left({ t}^{{ p}} +{ t}\right)\le \left({ 1}+{ t}\right)^{{ p}+{ 1}} {,}\quad { t}\ge { 0,}\quad { 0}<{ p}\le { 3}
$$
for ${ t}=\frac{x_i-x_i^{j_{i}}}{x_i^{j_{i}+1}-x_i}$ and
$1\le { p}\le { 3}$, we have
\begin{eqnarray*}
\alpha(x_i)\leq \frac{m_i}{2^p}(x_i^{j_{i}+1}-x_i)^{p+1}\left(1+\frac{x_i-x_i^{j_{i}}}{x_i^{j_{i}+1}-x_i}\right)^{p+1}\\
=\frac{m_i}{2^p}(x_i^{j_{i}+1}-x_i^{j_{i}})^{p+1}=\frac{m_i}{2^p}\left(\frac 1 {m_i} \right)^{p+1}=\frac{1}{(2m_i)^p}.
\end{eqnarray*}
\indent Since $\Omega(\gamma)$ is non-decreasing, we obtain
\begin{equation*}
{\cal E}_{{\bf m}}(f,{\bf x})  \leq  \Omega\left(\sqrt[p]{\sum_{i=1}^n \frac{1}{(2m_i)^p}  }  \right) = \Omega\left(\frac{1}{2}\sqrt[p]{\sum_{i=1}^n \frac{1}{(m_i)^p}  }  \right).
\end{equation*}
Since last inequality holds for any $f\in C_{D,p}(\Omega)$, we have
\begin{equation}\label{result2}
{E}_{{\bf m}}(C_{D,p}(\Omega)) \leq \Omega\left(\frac{1}{2}\sqrt[p]{\sum_{i=1}^n \frac{1}{(m_i)^p}  }  \right).
\end{equation}
Finally, we consider an extremal function $f_0$, defined as follows
\begin{equation*}
f_{0}({\bf x}):=\Omega\left(\sqrt[p]{\sum_{i=1}^{n}\left(\min\{x_i-x_i^{j_{i}}, x_i^{j_{i}+1}-x_i  \}  \right)^p  }  \right),
\end{equation*}
where ${\bf x}\in D_{\bf j}$  and $f_0({\bf x})=0$ for ${\bf x}\notin D_{\bf j}$. From the way $f_0$ is defined, it is easy to see that $f_0 \in C_{D,p}(\Omega)$.  Note that $f_{0}({\bf x}^{\bf j})=0$ for ${\bf j}=(j_1,...,j_n)$, $j_i=0,...,m_i$, $i=1,...,n$.  This, along with the linearity of $S_{\bf m}(f,{\bf x})$ on each element of the partition, implies that $S_{\bf m}(f,{\bf x})=0$ and
\begin{equation}\label{equal}
\|{\cal E}_{{\bf m}}(f_0)\|_C=\|f_0\|_C.
\end{equation}
\indent Using (\ref{equal}), we obtain
\begin{multline}\label{equality2}
E_{{\bf m}}(C_{D,p}(\Omega))\geq  \|{\cal E}_{{\bf m}}(f_0)\|_C
 =  \|f_0\|_C
 \geq  f_0\left(\frac{1}{2{\bf m}}  \right) =  \Omega\left( \frac{1}{2}\sqrt[p]{\sum_{i=1}^{n} \frac{1}{m_i^p} } \right).
\end{multline}
\indent Combining (\ref{result2}) and (\ref{equality2}), we have
\[
E_{{\bf m}}(C_{D,p}(\Omega))=\Omega\left( \frac{1}{2}\sqrt[p]{\sum_{i=1}^{n} \frac{1}{m_i^p} } \right). \;\;\;\;\;\;\;\;\;\;\;\;\;\;\;\; \square
\]
}

\section{Divided differences}

 In order to state the remaining results of this paper, we need to recall the definitions and some properties of divided differences.

Given function $f$, we define
\begin{equation}\label{divdifdef}
\delta(f;{\bf x, y}):=\frac{f({\bf x})-f({\bf y})}{\prod_{i\in G} (x_i-y_i) }, \qquad {\bf x, y}\in \mathbb{R}^n,
\end{equation}
where $G=\{i{\in\{1,...,n\}}:\;\; x_i\neq y_i\}$.

{ We remind the reader that for given ${\bf r}\in \{0,1\}^n$, the set $M$ is defined as ${M=\{i:\;r_i=1\}}$.}
In the remaining sections we will express $S_{\bf m}^{({\bf r})}(f;{\bf x})$ in terms of the divided difference of the function $f$ for the terms $x_i$, where $i\in M$. Recall that from (\ref{divspline}) we have
\begin{eqnarray*}
S_{\bf m}^{({\bf r})}(f;{\bf x}) &=& \left(\prod_{i\in M}m_i\right)F_M({\bf x})\\
&=& \left(\prod_{i\in M}m_i\right)\sum_{l_1=0}^{1}...\sum_{l_n=0 }^{1}(-1)^{\sum_{i\in M}l_i}\left(\prod_{i\notin M}H_{l_i,j_i}(x_i) \right) f(
{\bf x}^{{\bf j}+{\bf l}}
),
\end{eqnarray*}
where ${\bf x}^{{\bf j}+{\bf l}}=(x_1^{j_{1}+l_{1}},...,x_n^{j_{n}+l_{n}})$ and $j_i=0,...,m_i-1$, {${\bf l}=(l_1,...,l_n)\in\{0,1\}^n$} for ${i=1,...,n}$.
Fixing point $x_a$ such that $a\in M$, we have
\begin{eqnarray*}
 S_{\bf m}^{({\bf r})}(f;{\bf x})& = & m_a\left(\prod_{i\in M;\;\;i\neq a}m_i\right)\sum_{l_1=0}^{1}...\sum_{l_{a-1}=0 }^{1}\sum_{l_{a+1}=0 }^{1}...\sum_{l_n=0 }^{1}(-1)^{\sum_{i\in M;\;\;i\neq a}l_i}\\
& & \times \left(\prod_{i\notin M}H_{l_i,j_i}(x_i) \right)f(x_1^{j_{1}+l_{1}},...,x_a^{j_{a}},...,x_n^{j_{n}+l_{n}})\\
&+ & m_a\left(\prod_{i\in M;\;\;i\neq a}m_i\right)\sum_{l_1=0}^{1}...\sum_{l_{a-1}=0 }^{1}\sum_{l_{a+1}=0 }^{1}...\sum_{l_n=0 }^{1}(-1)^{\sum_{i\in M;\;\;i\neq a}l_i+1}\\
&  & \times \left(\prod_{i\notin M}H_{l_i,j_i}(x_i) \right)f(x_1^{j_{1}+l_{1}},...,x_a^{j_{a}+1},...,x_n^{j_{n}+l_{n}})\\
& = & m_a\left(\prod_{i\in M;\;\;i\neq a}m_i\right)\sum_{l_1=0}^{1}...\sum_{l_{a-1}=0 }^{1}\sum_{l_{a+1}=0 }^{1}...\sum_{l_n=0 }^{1}(-1)^{\sum_{i\in M;\;\;i\neq a}l_i+1}\\
& & \times \left(\prod_{i\notin M}H_{l_i,j_i}(x_i) \right)\left[f(x_1^{j_{1}+l_{1}},...,x_a^{j_{a}+1},...,x_n^{j_{n}+l_{n}}) \right.\\
& & \left.  - f(x_1^{j_{1}+l_{1}},...,x_a^{j_{a}},...,x_n^{j_{n}+l_{n}})\pm f(x_1^{j_{1}+l_{1}},...,x_{a},...,x_n^{j_{n}+l_{n}}) \right].\\
\end{eqnarray*}
\indent Setting
\begin{eqnarray*}
{\bf v}& = & (x_1^{j_{1}+l_{1}},...,x_a^{j_{a}-l_a+1},...,x_n^{j_{n}+l_{n}})\\
{\bf w}& = &(x_1^{j_{1}+l_{1}},...,x_a^{j_{a}},...,x_n^{j_{n}+l_{n}}) \\
{\bf y}& = &(x_1^{j_{1}+l_{1}},...,x_a^{j_{a}+1},...,x_n^{j_{n}+l_{n}}) \\
{\bf z}& = & (x_1^{j_{1}+l_{1}},...,x_{a},...,x_n^{j_{n}+l_{n}})
\end{eqnarray*}
and using the divided difference defined in (\ref{divdifdef}), we have
\begin{eqnarray*}
& &S_{\bf m}^{({\bf r})}(f;{\bf x})\\
& = &  m_a\left(\prod_{i\in M;\;\;i\neq a}m_i\right)\sum_{l_1=0}^{1}...\sum_{l_{a-1}=0 }^{1}\sum_{l_{a+1}=0 }^{1}...\sum_{l_n=0 }^{1}(-1)^{\sum_{i\in M;\;\;i\neq a}l_i+1}\\
& & \times\left(\prod_{i\notin M}H_{l_i,j_i}(x_i) \right) \left[(x_a^{j_{a}+1}-x_a)\delta\left( f;{\bf y,z} \right) + (x_a-x_{a,j_a}) \delta \left( f; {\bf z,w} \right)\right]\\
& = &\left(\prod_{i\in M;\;\;i\neq a}m_i\right)\sum_{l_1=0}^{1}\sum_{l_{2}=0 }^{1}...\sum_{l_n=0 }^{1}(-1)^{\sum_{i\in M;\;\;i\neq a}l_i+1}H_{l_a,j_a}(x_a)\\
& &\times\left(\prod_{i\notin M}H_{l_i,j_i}(x_i) \right) \delta \left(f;{\bf v,z}\right).
\end{eqnarray*}
\indent Repeating this process for each $i\in M$ {and using the fact that $|M|$ is even, we eventually obtain that the sign of each term is positive.} We
define the vectors ${\bf q, p} \in D_{{\bf j}}$ as follows
\begin{equation*}
{\bf q}= (q_1,q_2,...,q_n):\;\; q_i=\left\{
\begin{array}{lr}
x_i^{j_{i}-l_i+1}, & \;\;\; \mbox{if } i\in M\\
x_i^{j_{i}+l_{i}}, & \;\;\; \mbox{if } i\notin M
\end{array}
  \right.
\end{equation*}
\begin{equation*}
{\bf p}= (p_1,p_2,...,p_n):\;\; p_i=\left\{
\begin{array}{lr}
x_{i}, & \;\;\; \mbox{if } i\in M\\
x_i^{j_{i}+l_{i}}, & \;\;\; \mbox{if } i\notin M
\end{array}.
  \right.
\end{equation*}
\indent Therefore, we obtain
\begin{equation*}
S_{\bf m}^{({\bf r})}(f;{\bf x})=\sum_{l_1=0}^{1}...\sum_{l_n=0 }^{1}\left(\prod_{i=1}^{n}H_{l_i,j_i}(x_i) \right)\delta\left(f; {\bf q,p}\right).
\end{equation*}
\indent Combining it with $(\ref{divspline})$ {and using the triangle inequality}, we obtain
\begin{equation*}
{\cal E}_{{\bf m}}^{{\bf r}}(f;{\bf x})\leq\sum_{l_1=0}^{1}...\sum_{l_n=0 }^{1}\left(\prod_{i=1}^{n}H_{l_i,j_i}(x_i) \right)\left|f^{{\bf (r)}}({\bf x})-\delta\left(f; {\bf q,p}\right) \right|.
\end{equation*}

Using \cite{Mikeladze}, write the divided difference in integral form
\begin{multline}\label{divdifint}
{\cal E}_{{\bf m}}^{{\bf r}}\left( f;{\bf x} \right)\\
 \leq  \sum_{l_1=0}^{1}...\sum_{l_n=0 }^{1}\left(\prod_{i=1}^{n}H_{l_i,j_i}(x_i) \right)\int_{R''}\left|f^{{\bf (r)}}({\bf x})-f^{{\bf (r)}}({\bf x^*}) \right|\left( \prod_{i\in M} d\alpha_i \right),
\end{multline}
where
\[
R'':=[0,1]^{|M|}
\]
and
\[
{\bf x^*}=(x_1^*,x_2^*,...,x_n^*),\;\; \mbox{with } x_i^*=
\left\{ \begin{array}{lr}
x_i+\alpha_i(x_i^{j_i-l_i+1}-x_i),& \mbox{if } i\in M\\
x_i^{j_{i}+l_{i}},& \mbox{if } i\notin M
\end{array} \right. .
\]

\section{Approximation of derivatives of functions by derivatives of splines}

In this sections, we look at approximating a function's derivative by the
derivative of the linear spline that is constructed to interpolate the function itself.
The next two theorems correspond to the first two theorems if we only consider the original function.\\
\indent In the following result, we estimate the error of approximation on the class
$C_D^{{\bf r}}(\Omega)$, denoted by ${ E}_{{\bf m}}^{{\bf r}}(C_D^{{\bf r}}(\Omega))$.
\begin{theorem}\label{thm4}
{Let ${\bf r}\in \{0,1\}^n$ be given and $M:=\{i:\;r_i=1\}$. Let also an arbitrary function $\Omega(\boldsymbol\gamma)$  of MC-type, concave with respect to {$\gamma_i$}, where $i\notin M$, be given.}
Then for ${\bf m}\in \mathbb{N}^n$ with $m_i\geq 2$, $i=1,...,n$, the error of approximation
on the class $C_D^{{\bf r}}(\Omega)$ is
\[
E_{{\bf m}}^{{\bf r}}\left( C_D^{{\bf r}}(\Omega) \right)=\left(\prod_{i\in M} m_i \right)\int_{R}\Omega({\bf h})\left( \prod_{i\in M} d\gamma_i \right),
\]
where
\[
{\bf h}=(h_1,h_2,...,h_n),\qquad  {\mbox with}\qquad h_i=
\left\{ \begin{array}{lr}
\gamma_i,& \mbox{if } i\in M\\
\frac{1}{2m_i}, & \mbox{if } i\notin M
\end{array} \right.
\]
and
{$R=\prod\limits_{i\in M}[0,\frac{1}{m_i}]$}.
\end{theorem}
{\bf Proof}\\
\indent For any $\bf x$ from an arbitrary
$D_{{\bf j}}:=\prod_{i=1}^{n}[x_i^{j_{i}},x_i^{j_{i}+1}]$
and for given $f\in C_D^{{\bf r}}(\Omega)$ using the presentation of the
 divided difference in integral form \eqref{divdifint} and Definition \ref{defclassmul}, we have
\begin{equation}\label{Eomega}
\left| {\cal E}_{{\bf m}}^{{\bf r}}\left( f;{\bf x} \right) \right| \leq \sum_{l_1=0}^{1}...\sum_{l_n=0 }^{1}\left(\prod_{i=1}^{n}H_{l_i,j_i}(x_i) \right)\int_{R''}\Omega(\boldsymbol\beta)\left( \prod_{i\in M} d\alpha_i \right),
\end{equation}
where
\[
\boldsymbol\beta=(\beta_1, \beta_2,...,\beta_n),\mbox{ with } \beta_i= \left\{
\begin{array}{lr}
\alpha_i|x_i^{j_i-l_i+1}-x_i|,& \mbox{if } i\in M\\
|x_i-x_i^{j_{i}+l_{i}}|,& \mbox{if } i\notin M
\end{array} \right. .
\]
\indent By (\ref{Eomega}) and  for ${\bf x}\in D_{{\bf j}}$ we have
\begin{eqnarray}\label{R''}
\left| {\cal E}_{{\bf m}}^{{\bf r}}\left( f;{\bf x} \right) \right| & \leq & \sum_{l_1=0}^{1}...\sum_{l_n=0 }^{1}\left(\prod_{i=1}^{n}H_{l_i,j_i}(x_i) \right)\int_{R''}\Omega(\boldsymbol\beta^*)\left( \prod_{i\in M} d\alpha_i \right),
\end{eqnarray}
where
\[
\boldsymbol\beta^*=(\beta_1^*,\beta_2^*,...,\beta_n^*),\mbox{ with } \beta_i^*= \left\{
\begin{array}{lr}
\frac{\alpha_i}{m_i}H_{l_i,j_i}(x_i), & \mbox{if } i\in M\\
\frac{1}{m_i} H_{1-l_i,j_i}(x_i), & \mbox{if } i\notin M
\end{array} \right. .
\]
\indent Taking into account that $\Omega(\boldsymbol\gamma)$ is concave with respect to $\gamma_i$, $i\notin M$,
and performing a change of variables, we obtain
\begin{eqnarray*}
\left| {\cal E}_{{\bf m}}^{{\bf r}}\left( f;{\bf x} \right) \right| & \leq & \left(\prod_{i\in M} m_i \right) \sum_{l_{i_1}=0}^{1}...\sum_{l_{i_{|M|}}=0}^{1}\int_{R'}\Omega(\boldsymbol\beta^{'})\left( \prod_{i\in M} d\gamma_i \right),
\end{eqnarray*}
where {$i_k\in M$},
$R':=\prod\limits_{i\in M}[0,m_i^{-1}H_{l_i,j_i}(x_i)],$
 and
\[
\boldsymbol\beta^{'}=(\beta_1^{'},\beta_2^{'},...,\beta_n^{'}), \qquad \mbox{ with }\qquad \beta_i^{'}= \left\{
\begin{array}{lr}
\gamma_i, & \mbox{if } i\in M\\
\frac{2}{m_i} H_{0,j_i}(x_i)H_{1,j_i}(x_i), & \mbox{if } i\notin M
\end{array} \right. .
\]
\indent We have $\frac{2}{m_i} H_{0,j_i}(x_i)H_{1,j_i}(x_i)=\lambda(x_i)$, where $\lambda(x_i)$ is defined in (\ref{lambda}). Therefore, by (\ref{lambdamax}) we have that $\beta^{'}_i\leq \beta^{''}_i$ for all $i=1,...,n$, where
\begin{equation}\label{betageq}
\boldsymbol\beta^{''}=(\beta_1^{''},\beta_2^{''},...,\beta_n^{''}),\qquad \mbox{ with } \qquad \beta_i^{''}= \left\{
\begin{array}{lr}
\gamma_i, & \mbox{if } i\in M\\
\frac{1}{2m_i},  & \mbox{if } i\notin M
\end{array} \right. .
\end{equation}
\indent The fact that $\Omega(\boldsymbol\gamma)$ is non-decreasing together with (\ref{betageq}) imply
\begin{eqnarray}\label{b''}
\left| {\cal E}_{{\bf m}}^{{\bf r}}\left( f;{\bf x} \right) \right| & \leq & \left(\prod_{i\in M} m_i \right) \sum_{l_{i_1}=0}^{1}...\sum_{l_{i_{|M|}}=0}^{1}\int_{R'}\Omega(\boldsymbol\beta^{''})\left( \prod_{i\in M} d\gamma_i \right)\nonumber\\
& = & \Phi({\bf x}).
\end{eqnarray}
\indent It is easy to verify that $\Phi({\bf x})$ is continuous on $[x_i^{j_{i}},x_i^{j_{i}+1}]$ for $i\in M$. As $\Omega(\boldsymbol\gamma)\geq 0$, it is also easy to verify that
\begin{equation}\label{Phi}
\max\{\Phi({\bf x}):\; x_i\in[x_i^{j_{i}},x_i^{j_{i}+1}],i\in M \} =\Phi({\bf x}^{{\bf j}+{\bf l}}),{\bf l}=(l_1,...,l_n), l_i=0,1.
\end{equation}
\indent By (\ref{b''}) and (\ref{Phi}), we have
\begin{eqnarray*}
\left| {\cal E}_{{\bf m}}^{{\bf r}}\left( f;{\bf x}\right) \right| & \leq & \Phi({\bf x}^{\bf j})
\nonumber\\
& = &  \left(\prod_{i\in M} m_i \right) 
\int_{R}\Omega(\boldsymbol\beta^{''})\left( \prod_{i\in M} d\gamma_i \right),
\end{eqnarray*}
where $R=\prod\limits_{i\in M}[0,\frac{1}{m_i}]$. Setting $\beta^{''}={\bf h}$ {and using the fact that the last inequality holds for any function $f\in C_D^{{\bf r}}(\Omega)$}, we have the upper bound
\begin{equation}\label{result4}
E_{{\bf m}}^{{\bf r}}\left( C_D^{{\bf r}}(\Omega) \right) \leq \left(\prod_{i\in M} m_i \right)\int_{R}\Omega({\bf h})\left( \prod_{i\in M} d\gamma_i \right).
\end{equation}
\indent In order to show that the upper bound (\ref{result4}) is achieved, we introduce the extremal function
\begin{equation*}
d({\bf x}):=\int_S\pi({\bf h^{'}})\left(\prod_{i\in M}d \gamma_i \right),
\end{equation*}
where $S=\prod\limits_{i\in M}[0,x_i]$,
 \[
{\bf h^{'}}=(h^{'}_1, h^{'}_2,..., h^{'}_n),\qquad \mbox{with }  \qquad h^{'}_i= \left\{
\begin{array}{lr}
\gamma_i, & \mbox{if } i\in M\\
x_i,  & \mbox{if } i\notin M
\end{array} \right. ,
\]
for ${\bf{x}\in}\prod\limits_{i\in M}[\frac {l_i} {m_i},\frac {l_i+1} {m_i}]\times\prod\limits_{i\notin M}[\frac{l_i}{2m_i},\frac {l_i+1} {2m_i}]$, $l_i=0,1$, $i=1,...,n$,
\begin{equation*}
\pi({\bf x})=\Omega({\bf h}^{''})-\left(\prod_{i\in M} m_i \right)\int_{R}\Omega({\bf h})\left( \prod_{i\in M} d\gamma_i \right),
\end{equation*}
where
 \[
{\bf h^{''}}=(h^{''}_1,h^{''}_2,..., h^{''}_n),\qquad \mbox{with }  \qquad h^{''}_i= \left\{
\begin{array}{lr}
(-1)^{l_i}\left(\frac{1}{m_i}-x_i\right), & \mbox{if } i\in M\\
{(-1)^{l_i}\left(\frac{1}{2m_i}-x_i\right),}
& \mbox{if } i\notin M
\end{array} \right.
\]
{and extend function $\pi({\bf x})$ to be $\frac{2}{m_i}$-periodic for $x_i$ such that $i\in M$, and to be $\frac{1}{m_i}$-periodic for $x_i$ where $i\notin M$.}

From the way that $d({\bf x})$ is defined, we see that it is in the class $C_D^{{\bf r}}(\Omega)$.  In addition, we have that $d({\bf x}^{\bf j})=0$ for ${\bf j}=(j_1,...,j_n)$, $j_i=0,...,m_i$, so it follows that $S_{\bf m}(d;{\bf x})=0$ as the spline is linear in each $x_i$ on every element of the partition.  With $S^{{(\bf r)}}_{\bf m}(d;{\bf x})=0$, we obtain the following inequality
\begin{eqnarray}\label{lower4}
E_{{\bf m}}^{{\bf r}}\left( C_D^{{\bf r}}(\Omega)\right)& \geq & \|{\cal E}_{{\bf m}}^{{\bf r}}(d)\|_C
 =  \|d^{({\bf r})}
\|_C
 \geq  \left| d^{({\bf r})}\left( \boldsymbol\phi \right) \right| \nonumber \\
& = & \left(\prod_{i\in M} m_i \right)\int_{R}\Omega({\bf h})\left( \prod_{i\in M} d\gamma_i \right),
\end{eqnarray}
where
\[
\boldsymbol\phi=(\phi_1,\phi_2,..., \phi_n),\qquad \mbox{ with }\qquad  \phi_i= \left\{
\begin{array}{lr}
\frac{1}{m_i}, & \mbox{if } i\in M\\
\frac{1}{2m_i},  & \mbox{if } i\notin M
\end{array} \right. .
\]
\indent Comparing (\ref{result4}) and (\ref{lower4}), we obtain the approximation error to be
\[
E_{{\bf m}}^{{\bf r}}\left( C_D^{{\bf r}}(\Omega) \right)=\left(\prod_{i\in M} m_i \right)\int_{R}\Omega({\bf h})\left( \prod_{i\in M} d\gamma_i \right).\;\;\;\;\;\;\;\;\;\square
\]
\indent Theorem \ref{thm5} provides the estimate for the error of approximation on the class
$C_{D,p}^{{\bf r}}(\Omega)$.
\begin{theorem}\label{thm5}
Let ${\bf r}\in \{0,1\}^n$ be given and let $M=\{i:\;r_i=1\}$. Let also an arbitrary univariate, concave, MC-type function $\Omega(\gamma)$ be given.  Then for ${\bf m}\in \mathbb{N}^n$ with $m_i\geq 2$, $i=1,...,n$, the error of approximation on the class $C_{D,p}^{{\bf r}}(\Omega)$, $1\le p\le3$, is
\small{\begin{eqnarray*}
E_{{\bf m}}^{{\bf r}}\left( C_{D,p}^{{\bf r}}(\Omega) \right)= \left( \prod_{i\in M}m_i \right)\int_R\Omega\left(\sqrt[p]{\left[\sum_{i\in M} \gamma_i^p \right]+\left[\sum_{i\notin M} \frac{1}{(2m_i)^p}  \right] } \right)\left( \prod_{i\in M}d \gamma_i  \right),
\end{eqnarray*}}
where $R=\prod\limits_{i\in M}[0,\frac{1}{m_i}]$.
\end{theorem}
{\bf Proof}\\
Let arbitrary $f\in C_D^{{\bf r}}(\Omega)$ be given.  For any $\bf x$ from an arbitrary
$D_{{\bf j}}:=\prod_{i=1}^{n}[x_i^{j_{i}},x_i^{j_{i}+1}]$,
using the estimate of the error in the form (\ref{divdifint}), and using Definitions \ref{defombig} of $\Omega(\gamma)$ and Definition \ref{defclassone} of the class $C_D^{{\bf r}}(\Omega)$,
we have
\begin{multline*}
\left| {\cal E}_{{\bf m}}^{{\bf r}}\left( f;{\bf x} \right) \right|  \leq  \sum_{l_1=0}^{1}...\sum_{l_n=0 }^{1}\left(\prod_{i=1}^{n}H_{l_i,j_i}(x_i) \right)\\
\times \int_{R''} \omega_{p}\left(f^{({\bf r})};\sqrt[p]{ \sum_{i\in M}\left[ \alpha_i|x_i^{j_i-l_i+1}-x_i|\right]^p+\sum_{i\notin M} |x_i-x_i^{j_{i}+l_{i}}|^p } \right)\left( \prod_{i\in M} d\alpha_i \right) \\
 \leq  \sum_{l_1=0}^{1}...\sum_{l_n=0 }^{1}\left(\prod_{i=1}^{n}H_{l_i,j_i}(x_i) \right) \\
\times \int_{R''}\Omega\left(\sqrt[p]{ \sum_{i\in M}\left( \frac{\alpha_i}{m_i}H_{l_i,j_i}\right)^p+\sum_{i\notin M} \left(\frac 1 {m_i}H_{1-l_i,j_i}\right)^p } \right)\left( \prod_{i\in M} d\alpha_i \right),\\
\end{multline*}
where $R'':=[0,1]^{|M|}.$

\indent  Performing the change of variables, taking into account that $\Omega(\gamma)$ is concave,
and using notation $R':=\prod\limits_{i\in M}[0,m_i^{-1}H_{l_i,j_i}(x_i)]$
, we have
\begin{eqnarray*}
\left| {\cal E}_{{\bf m}}^{{\bf r}}\left( f;{\bf x} \right) \right| & \leq & \left(\prod_{i\in M} m_i \right) \sum_{l_1=0}^{1}...\sum_{l_n=0 }^{1}\left(\prod_{i\notin M}H_{l_i,j_i}(x_i) \right) \nonumber\\
& & \times \int_{R'}\Omega\left(\sqrt[p]{ \left[\sum_{i\in M} \gamma_i^p\right]+\left[\sum_{i\notin M} \left(\frac 1 {m_i}H_{1-l_i,j_i}\right)^p \right] } \right)\left( \prod_{i\in M} d\gamma_i \right)\nonumber\\
&\leq & \left(\prod_{i\in M} m_i \right) \sum_{l_{i_1}=0}^{1}...\sum_{l_{i_{|M|}}=0}^{1}\\
& &\times\int_{R'} \Omega\left(\sqrt[p]{ \left[\sum_{i\in M} \gamma_i^p\right]+\left[\sum_{i\notin M} \alpha^p(x_i) \right] } \right)\left( \prod_{i\in M} d\gamma_i \right),
\end{eqnarray*}
where {$i_k\in M$}
and $\alpha(x_i)$ is defined in (\ref{alpha}).\\
\indent Recall that in Theorem \ref{thm2} we have proved that for any ${\bf x}\in D_{\bf j}$
$$
\alpha(x_i)\leq \frac{1}{(2m_i)^p},\;i=1,...,n.
$$
%
\indent Therefore, since $\Omega(\gamma)$ is non-decreasing, we have
\begin{eqnarray}\label{premax}
\left| {\cal E}_{{\bf m}}^{{\bf r}}\left( f;{\bf x} \right) \right| & \leq & \left(\prod_{i\in M} m_i \right)
\sum_{l_{i_1}=0}^{1}...\sum_{l_{i_{|M|}}=0}^{1}\\
& & \times \int_{R'}\Omega\left(\sqrt[p]{\left[\sum_{i\in M} \gamma_i^p\right]+\left[\sum_{i\notin M} \frac{1}{(2m_i)^p} \right] } \right)\left( \prod_{i\in M} d\gamma_i \right)\nonumber\\
& = & \mu({\bf x}). \nonumber
\end{eqnarray}
\indent It is easy to verify that
\begin{multline}\label{maxmu}
\max\{ \mu({\bf x}): \;\;{\bf x}\in D_{{\bf j}} \}\\
=\left(\prod_{i\in M} m_i \right)\int_{R}\Omega\left(\sqrt[p]{\left[\sum_{i\in M} \gamma_i^p\right]+\left[\sum_{i\notin M} \frac{1}{(2m_i)^p} \right] } \right)\left( \prod_{i\in M} d\gamma_i \right),
\end{multline}
where $R=[0,\frac{1}{m_i}]^{|M|}$ for $i\in M$. Therefore, using (\ref{maxmu}) and the fact that (\ref{premax}) holds true for any $f\in  C_{D,p}^{{\bf r}}(\Omega)$, we have
\begin{multline}\label{result5}
E_{{\bf m}}^{{\bf r}}\left( C_{D,p}^{{\bf r}}(\Omega) \right)
\leq\left( \prod_{i\in M}m_i \right)\int_R\Omega\left(\sqrt[p]{\sum_{i\in M} \gamma_i^p+\sum_{i\notin M} \frac{1}{(2m_i)^p} } \right)\left( \prod_{i\in M}d \gamma_i  \right).
\end{multline}
\indent In order to show that equality in (\ref{result5}) is achieved, we introduce the extremal function
\begin{equation*}
e({\bf x}):=\int_S\pi({\bf h^{'}})\left(\prod_{i\in M}d \gamma_i \right)
\end{equation*}
where $S=[0,x_i]^{|M|}$ for $i\in M$,
 \[
{\bf h^{'}}=(h^{'}_1, h^{'}_2,..., h^{'}_n),\qquad \mbox{with }  \qquad h^{'}_i= \left\{
\begin{array}{lr}
\gamma_i, & \mbox{if } i\in M\\
x_i,  & \mbox{if } i\notin M
\end{array} \right. ,
\]
{and function $\pi(\bf x)$
is defined as
\begin{multline*}
\pi({\bf x}) \\
= \Omega\left(\sqrt[p]{ \left[\sum_{i\in M}\left( (-1)^{l_i}\left(\frac{1}{m_i}-x_i \right)
\right)^p\right]+ \left[\sum_{i\notin M} \left((-1)^{l_i}\left(\frac{1}{2m_i}-x_i\right)\right)^p
\right]} \right)\\
 -\left( \prod_{i\in M}m_i \right)\int_R\Omega\left(\sqrt[p]{\left[\sum_{i\in M} \gamma_i^p\right]+\left[\sum_{i\notin M} \frac{1}{(2m_i)^p} \right] } \right)\left( \prod_{i\in M}d \gamma_i  \right)
\end{multline*}
for ${\bf{x}\in}\prod\limits_{i\in M}[\frac {l_i} {m_i},\frac {l_i+1} {m_i}]\times\prod\limits_{i\notin M}[\frac{l_i}{2m_i},\frac {l_i+1} {2m_i}]$, $l_i=0,1$, $i=1,...,n$ and then extended so that $\pi({\bf x})$ is $\frac{2}{m_i}$-periodic for $x_i$ such that $i\in M$ and is $\frac{1}{m_i}$-periodic for $x_i$ where $i\notin M$.}

\indent From the way $e({\bf x})$ is defined, we see that it clearly belongs to the class $C_D^{{\bf r}}(\Omega)$.  We have that $e({\bf x}^{\bf j})=0$ for ${\bf j}=(j_1,...,j_n)$, $j_i=0,...,m_i$, so it follows that $S_{\bf m}(e;{\bf x})=0$ as the spline is linear in each $x_i$ on each partition element.  With $S^{({\bf r})}_{\bf m}(e;{\bf x})=0$, we obtain the following inequality
\begin{multline}\label{lower5}
E_{{\bf m}}^{{\bf r}}\left( C_{D,p}^{{\bf r}}(\Omega)\right) \geq  \|{\cal E}_{{\bf m}}^{{\bf r}}(e)\|_C = \|e^{({\bf r})}\|_C \geq \left| e^{({\bf r})}\left( \boldsymbol\phi \right) \right| \\
 =  \left( \prod_{i\in M}m_i \right)\int_R\Omega\left(\sqrt[p]{\left[\sum_{i\in M} \gamma_i^p\right]+\left[\sum_{i\notin M} \frac{1}{(2m_i)^p} \right] } \right)\left( \prod_{i\in M}d \gamma_i  \right),
\end{multline}
where
\[
\boldsymbol\phi=(\phi_1,\phi_2,..., \phi_n),\qquad \mbox{with } \;\; \phi_i= \left\{
\begin{array}{lr}
\frac{1}{m_i}, & \mbox{if } i\in M\\
\frac{1}{2m_i},  & \mbox{if } i\notin M
\end{array} \right. .
\]
\indent Comparing (\ref{result5}) and (\ref{lower5}), we obtain the following error of approximation
\begin{multline*}
E_{{\bf m}}^{{\bf r}}\left( C_{D,\rho_p}^{{\bf r}}(\Omega) \right)\\
=\left( \prod_{i\in M}m_i \right)\int_R\Omega\left(\sqrt[p]{\left[\sum_{i\in M} \gamma_i^p\right]+\left[\sum_{i\notin M} \frac{1}{(2m_i)^p} \right] } \right)\left( \prod_{i\in M}d \gamma_i  \right). \;\;\;\square
\end{multline*}

{\bf Remark.} If in the statement of the theorem all coordinates of the vector ${\bf r}$ are equal to $1$, then the assumption on $\Omega(\gamma)$ to be concave can be removed.


\begin{thebibliography}{99}
\bibitem{Malozemov66}
V.\,N.~Malozemov, On the deviation of broken lines, [in Russian] {\it Vestn. Leningr. Univ. Math. Mekh. Astron}, {\bf 7} (1966), pp. 150--153.

\bibitem{Malozemov67}
V.\,N.~Malozemov, On polygonal interpolation, {\it Mathematical Notes},
{\bf 1}:5 (1967), pp. 537--540.

\bibitem{Storchai75}
V.\,F.~Storchai, Approximation of continuous functions of two variables by
polyhedral functions and splines in the unform metric, in: {\it Investigations in Contemporary Problems of Summation and Approximation of Functions and Their Applications} [in Russian], Dnepropetrovsk State University, Dnepropetrovsk (1975), pp. 82--89.

\bibitem{Vakarchuk90}
S.\,B.~Vakarchuk, Interpolation by bilinear splines, {\it Mathematical Notes}, {\bf 47}:5 (1990), pp. 441--444.

\bibitem{Shabozov}
M.\,Sh.~Shabozov, On the error of the interpolation by bilinear splines, {\it Ukrainian Mathematical Journal}, {\bf 46}:11 (1994), pp. 1719--1726.

\bibitem{Vakarchuk05}
S.\,B.~Vakarchuk, K.\,Yu.~Myskin, Some problems of simultaneous approximation of functions of two variables and their derivatives by interpolation bilinear splines, {\it Ukrainian Mathematical Journal}, {\bf 57}:2 (2005), pp. 173--185.

\bibitem{Storchai72}
V.\,F.~Storchai, Approximation of continuous functions of two variables by
splines in metric $C$, in: {\it Investigations in Contemporary Problems of Summation and Approximation of Functions and Their Applications} [in Russian], Dnepropetrovsk State University, Dnepropetrovsk (1972), pp. 66--68.

\bibitem{Korneychuk}
N.\,P.~Korneychuk, {\it Sharp constants in approximation theory}, Encyclopedia of Mathematics and its Applications V. 38, Cambridge University Press, 1991.

\bibitem{Mikeladze}
Sh.\,E.~Mikeladze, {\it Numerical Methods of Mathematical Analysis } [in Russian], Gostekhizdat, Moscow, 1953.









\end{thebibliography}
\end{document}